\documentclass[12pt]{amsart}
\usepackage{epsfig}
\usepackage{graphics}
\numberwithin{equation}{section}
\theoremstyle{plain}
   \newtheorem{thm}{Theorem}[section]

   \newtheorem{lemma}[thm]{Lemma}
   \newtheorem{claim}[thm]{Claim}

   \newtheorem{conj}[thm]{Conjecture}

\theoremstyle{definition}
   \newtheorem{dfn}[thm]{Definition}
   
\theoremstyle{remark}

\newcommand{\C}{\mathbb C}

\newcommand{\R}{\mathbb R}

\newcommand{\bd}{\partial}

\begin{document}
\title[Positive open book decompositions along arbitrary links]{Stein fillable $3$-manifolds admit positive open book decompositions along arbitrary links}
\author{Masaharu Ishikawa}
\address{Department of Mathematics, Tokyo Institute of Technology,\\
2-12-1, Oh-okayama, Meguro-ku, Tokyo, 152-8551, Japan}
\email{ishikawa@math.titech.ac.jp}

\begin{abstract}
It is known by A.~Loi and R.~Piergallini that a closed, oriented, smooth
$3$-manifold is Stein fillable if and only if it has
a positive open book decomposition. In the present paper
we will show that for every link $L$ in a Stein fillable $3$-manifold
there exists an additional knot $L'$ to $L$ such that the link
$L\cup L'$ is the binding of a positive open book decomposition of
the Stein fillable $3$-manifold. To prove the assertion,
we will use the divide, which is a generalization
of real morsification theory of complex plane curve singularities,
and $2$-handle attachings along Legendrian curves.
\end{abstract}

\maketitle

\section{Introduction}

An open book decomposition of a closed, oriented, smooth $3$-manifold
is the following:
Let $S$ be a compact, oriented, smooth $2$-dimensional manifold
with boundary $\bd S$
and $h$ an automorphism of $S$ which is the identity on $\bd S$.
If a closed, oriented, smooth $3$-manifold $M$ can be obtained from
$S\times [0,1]$ by identifying $(h(x),0)$ and $(x,1)$ for $x\in S$
and $(y,0)$ and $(y,t)$ for $y\in\bd S$ and all $t\in [0,1]$,
then we say that $M$ has an {\it open book decomposition}.
The manifold $S$ is called a {\it fiber surface} and $\bd S$ the {\it binding}.
Note that the binding is a {\it fibered link} in $M$.
An open book decomposition of $M$ is said {\it positive}
if its monodromy $h$ consists of a product of positive Dehn twists.

In~\cite{stallings:78} J.R.~Stallings proved that for every link $L$
in $S^3$ there exists a knot $L'$ such that $L\cup L'$
is the binding of an open book decomposition of $S^3$.
Moreover he proved that the knot $L'$ can be chosen
in such a way that $L'$ and each connected component of $L$
has an arbitrarily prescribed linking number.
As a corollary, we can conclude that for every link $L$
in a closed, oriented, smooth $3$-manifold $M$
there exists a knot $L'$ such that $L\cup L'$
is the binding of an open book decomposition of $M$,
see~\cite[Theorem~7.4]{i}.

Now we focus on positive open book decompositions 
of closed, oriented, smooth $3$-manifolds.
The existence of a positive open book decomposition is related to
the following Stein fillability.
A complex manifold $X$ is called {\it Stein}
if it admits a strictly plurisubharmonic function $\phi:X\to\R$
which is proper and bounded below.
Each level set $\phi^{-1}(c)$, for $c\in\R$, is strictly pseudoconvex,
where $\phi^{-1}(c)$ is oriented as the boundary of $\phi^{-1}((-\infty,c])$.
In the case where the complex dimension of $X$ is $2$,
we call $\phi^{-1}((-\infty,c])$ a {\it compact Stein surface with boundary},
or we may call it simply a {\it compact Stein surface}.
A closed, oriented, smooth $3$-manifold $M$ is called
{\it Stein fillable} if there exists a compact Stein surface
$\phi^{-1}((-\infty,c])$ such that $M=\phi^{-1}(c)$.
In the paper~\cite{eliashberg:90}, Y.~Eliashberg gave a topological
characterization of compact Stein surfaces
using handle decompositions along Legendrian curves.
Then R.E.~Gompf described, in~\cite{gompf:98}, the handle
decompositions explicitly using framed Legendrian links and Kirby calculus.
Following these studies, A.~Loi and R.~Piergallini represented
the manifold $\phi^{-1}((-\infty,c])$ as branched covers of $B^4$
using framed Legendrian link presentations of handle decompositions.
As a consequence, they proved that a closed, oriented, smooth $3$-manifold
is Stein fillable if and only if it has
a positive open book decomposition\cite{lp:2001}.
Note that a very few examples of non-Stein fillable $3$-manifolds
are known, see~\cite{lisca:98,lisca:99}.

In the present paper we study what kinds of positive open book
decompositions a closed, oriented, smooth $3$-manifold has.
Because of the result of Loi and Piergallini, we can exclude
non-Stein fillable $3$-manifolds.
The main result in this paper is the following:

\begin{thm}\label{mainthm}
Let $M$ be a Stein fillable $3$-manifold and fix a link $L$ in $M$. Then
there exists a knot $L'$ in $M\setminus L$ such that
the union $L\cup L'$ is the binding of
a positive open book decomposition of $M$.
\end{thm}

In~\cite{gi}, W.~Gibson and the author proved the theorem
for the case where $M$ is $S^3$, and then in~\cite{i}
the author proved it for the case where
$M$ is either $\#nS^1\times S^2$ or the unit tangent bundles
to closed, oriented surfaces.
Since a non-Stein fillable $3$-manifold does not admit any positive
open book decomposition, the result in the present paper
is the final solution of this series of studies.

To prove the main theorem we will use {\it divides} as we did
in~\cite{gi} and~\cite{i}. A divide is the image of a generic,
relative immersion of a finite number of copies of the unit interval
or the unit circle into the genus-$g$ surface
$\Sigma_{g,k}$ with $k$ boundary components.
It was originally defined in the unit disk by N.~A'Campo in~\cite{acampo:99}
as a generalization of real morsification theory of complex plane curve
singularities\cite{acampo:75,acampo:74,gz1,gz2,gz3}.
Each divide in the unit disk determines a link in $S^3$.
In~\cite{acampo:99} A'Campo proved that this link is fibered
with positive monodromy diffeomorphism if the divide is connected.
To each divide in $\Sigma_{g,k}$ is assigned a link in
the unit tangent bundle to $\Sigma_{g,k}$ if $k=0$,
$S^3$ if $2g+k-1=0$, and the connected sum of $2g+k-1$ copies of
$S^1\times S^2$ if $2g+k-1>0$.
In~\cite{i}, it is proved that the link is fibered
if the divide satisfies a certain condition.
In the present paper, we deal with only divides in $\Sigma_{0,n+1}$,
whose links are defined in $S^3$ if $n=0$ and in $\#nS^1\times S^2$ if $n>0$.
Hereafter we set $\#0S^1\times S^2=S^3$ by definition.
The fibration theorem in~\cite{i} is proved by constructing
a Lefschetz fibration from a compact, connected,
oriented $4$-manifold to a disk such that the fibration of a divide
is induced on its boundary.

For constructing the expected positive open book decomposition of a
Stein fillable $3$-manifold
we first prepare, using a divide in $\Sigma_{0,n+1}$,
a Lefschetz fibration over a disk which induces, on its boundary,
a positive open book decomposition of $\#nS^1\times S^2$
whose binding contains the prescribed link.
Due to the work of Eliashberg in~\cite{eliashberg:90},
it is known that a compact Stein surface
can be obtained from a $4$-ball with $1$-handles by 
attaching $2$-handles along Legendrian knots with coefficients
the canonical framing minus $1$ (see~\cite{gompf:98}).
Note that the boundary of the $4$-ball with $1$-handles is $\#nS^1\times S^2$.
For yielding a Lefschetz fibration after the $2$-handle attachings again,
we use the method of S.~Akbulut and B.~Ozbagci in~\cite{ao:2001},
that is, we will find a good position of a fiber surface of
the positive open book decomposition appearing on the boundary
of a Lefschetz fibration
such that (i) the links for $2$-handle attachings lie on the fiber surface and
(ii) their canonical framings coincide with
the product framings of the Lefschetz fibration.
These two conditions ensure that the $2$-handle attachings
yield a Lefschetz fibration again
and thus, on its boundary, we obtain the Stein fillable $3$-manifold
equipped with the expected positive open book decomposition.

This paper is organized as follows:
In Section~2 we introduce the unit disk $D_n$ in $\R^2$ with $n$ holes
and define the $3$-manifold $\#nS^1\times S^2$ as the boundary of
a small compact neighborhood of $D_n$ in $\C^2$.
In Section~3 we introduce divides in $D_n$ and associated
positive open book decompositions of $\#nS^1\times S^2$.
In the end of this section we briefly explain how to see
the fiber surfaces of the positive open book decompositions associated
with divides.
In Section~4 we shortly introduce oriented divides.
Section~5 is devoted to the proof of the main theorem.
In Section~6 we propose a conjecture for the future.

I would like to thank Professor Norbert A'Campo for many useful
conversations and suggestions in the preliminary stage of this study.

\section{Wave fronts on the unit disk with $n$ holes}

Let $(x_1+iu_1, x_2+iu_2)\in\C^2$ be the complex valued coordinates of $\C^2$
and set $D:=\{(x_1,x_2)\in\R^2 \mid x_1^2+x_2^2\leq 1\}\subset\R^2\subset\C^2$.
We denote by $D_n$ the unit disk with $n$ holes obtained
from $D$ by removing $n$ open disks, where $n\geq 0$, and
by $\bd D_n$ the boundary of $D_n$.

We thicken $D_n$ in $\C^2$ as follows:
Let $A_\delta$ be a small compact neighborhood of $\bd D_n$ in $D_n$
with width $\delta>0$ and suppose that $\delta$ is sufficiently small.
Set $B_\delta:=D_n\setminus A_\delta$ and thicken it as
$
   \hat B_\delta:=\{x+iu\in\C^2\mid x\in B_\delta,\, |u|\leq\delta\}
$.
Set $\alpha_\delta:=\bd A_\delta\setminus\bd D_n$
and thicken $A_\delta$ as
\[
   \hat A_\delta:=\{x+iu\in\C^2\mid
      x\in A_\delta,\, |u|^2+d(x,\alpha_\delta)^2\leq \delta^2\},
\]
where $d(\cdot,\cdot)$ is the minimal distance of the sets.
We then define a thickened disk $N_\delta(D_n)$
of $D_n$ in $\C^2$ by the union
$N_\delta(D_n):=\hat A_\delta\cup\hat B_\delta$,
which is a compact $4$-dimensional manifold in $\C^2$.
We denote the boundary of $N_\delta(D_n)$ by $M_n$.
It is easy to see that $M_n$ is $\#nS^1\times S^2$,
which is the connected sum of $n$ copies of $S^1\times S^2$ if $n\geq 1$
and $S^3$ if $n=0$.

According to the decomposition $N_\delta(D_n)=\hat A_\delta\cup\hat B_\delta$,
we decompose $M_n$ into two pieces, $\bd N_A:=\bd\hat A_\delta\cap M_n$
and $\bd N_B:=\bd\hat B_\delta\cap M_n$. The first piece $\bd N_A$
is homeomorphic to a union of solid tori whose core curves are $\bd D_n$
and the second piece $\bd N_B$ is $B_\delta\times S^1$.
We may regard the second piece as the unit tangent bundle to $B_\delta$.
For the sake of simplicity, we will abbreviate $B_\delta$ to $B$.

When discussing a contact structure in $M_n$
we use the $3$-manifold $M_n':=\bd N_A'\cup\bd N_B$ instead of $M_n$,
where $\bd N_A':=\{x+iu\in M_n\mid
x\in A_\delta,\, |u|^2+d(x,\alpha_\delta)^3=\delta^2\}$,
since $|u|=\delta$ and $|u|^2+d(x,\alpha_\delta)^2=\delta^2$ 
are not $C^2$-differentiable at each point in the boundary of $\bd N_A$.
The complex tangency of $M_n'\subset\C^2$ induces
a contact structure $\xi$ in $M_n'$ and we can check by direct calculation
that if $\delta>0$ is sufficiently small then $M_n'$ is strictly pseudoconvex.
Thus $\xi$ is the standard contact structure in $\#nS^1\times S^2$
compatible with its standard orientation.
Here the standard contact structure is the unique tight contact structure in
$\#nS^1\times S^2$ (up to isotopy).
Note that the uniqueness is known in~\cite{bennequin:83,eliashberg:92}.

A {\it contact element} on $B$ is a line tangent to $B$ at a point in $B$.
We consider co-oriented contact elements on $B$.
For each point $x\in B$ and local coordinates $(x_1,x_2)$ in $B$
centered at $x$, a co-oriented contact element is represented
by the equation
\[
   \alpha_x:=(cos\theta)dx_1+(\sin\theta)dx_2=0,
\]
where the direction $\theta$ is the co-orientation of the contact element.
The co-oriented contact element is a point in the unit cotangent bundle
to $B$ and it naturally corresponds to a unit tangent vector in the
unit tangent bundle to $B$. Thus the space of co-oriented contact elements
on $B$ corresponds to the unit tangent bundle to $B$, which
is the $3$-manifold $\bd N_B$.
In $\bd N_B$, the contact form $\alpha$ of $\xi$ is given by
$\alpha|_{\bd N_B}=-(u_1dx_1+u_2dx_2)$, where $u_1^2+u_2^2=\delta^2$,
and hence the contact structure $\xi$ of $\bd N_B\subset M_n'$
coincides with the contact structure
determined by the co-oriented contact elements on $B$.

A link in a $3$-manifold equipped with a contact structure
is called {\it Legendrian} if all the tangent vectors to the link
lie on the $2$-plane field of the contact structure.

A {\it wave front $w$} in $B$ is the image of a (generic) immersion 
of a finite number of copies of the unit circle into $B$,
possibly with cusps, equipped with co-orientation.
Each line tangent to $w$ with co-orientation is
a co-oriented contact element and those contact elements constitute
a Legendrian link $L$ in $\bd N_B$.
We say that $L$ is the {\it Legendrian link of $w$}, or
$w$ is a {\it wave front representative of $L$}.

\begin{lemma}\label{lemma1}
Every Legendrian link in $M_n$ has a wave front representative on $B$
up to Legendrian isotopy.
\end{lemma}

\begin{proof}
Since the piece $\bd N_A$ is the union of solid tori,
using Legendrian isotopy we can assume that the Legendrian link stays
in $\bd N_B$. Then the assertion is a well-known fact in contact topology.
\end{proof}

\section{Fibration theorem of divides in $D_n$}

In this section we introduce divides in $D_n$ and their links
and fibrations in $\#nS^1\times S^2$.
The idea is based on the works of A'Campo in~\cite{acampo:99,acampo:98a}.
Divides in the genus-$g$ surface, possibly with boundary, have been
studied in~\cite{i}, which contain divides in $D_n$.

\begin{dfn}\label{dfndiv}
A {\it divide} $P$ in $D_n$ is the image of a generic,
relative immersion of a finite number of copies of the unit interval or
the unit circle into $D_n$.
Each image of the unit interval (resp. circle) is called
an {\it interval (resp. circle) component} of $P$. The generic and relative
conditions are the following:
\begin{itemize}
   \item[(i)] the image has neither self-tangent points nor triple points;
   \item[(ii)] each endpoint of an interval component lies on
$\bd D_n$ and the interval component intersects $\bd D_n$
at the endpoints transversely;
   \item[(iii)] a circle component does not intersect $\bd D_n$.
\end{itemize}
\end{dfn}

An {\it edge} of $P$ is the closure of a connected component of
$P\setminus\{double ~points\}$ and a {\it region} of $P$ is a connected
component of $D_n\setminus P$. If a region of $P$ is bounded by
only $P$ then it is called an {\it interior region},
and otherwise it is called an {\it exterior region}.
For each exterior region, its intersection with $\bd D_n$
is called the {\it outside boundary}.

\begin{dfn}
A divide $P$ in $D_n$ is {\it admissible} if it
satisfies the following:
\begin{itemize}
   \item[(iv)] $P$ is connected;
   \item[(v)] each interior region of $P$ is simply connected;
   \item[(vi)] each exterior region of $P$ is either (a) simply connected
and the outside boundary in that region is connected\footnote{
   The condition `the outside boundary in that region is connected',
   which was not mentioned in~\cite{i},
   is necessary for drawing the level sets of
   an associated Morse function $f_P$.},
or (b)
an annulus such that one boundary component is a component of
$\bd\Sigma_{g,n}$ and the other is contained in $P$;
   \item[(vii)] $P$ allows a checkerboard coloring, that is, a coloring
with two colors, black and white, such that if two regions of $P$
share an edge of $P$ in their boundaries then they are painted with
different colors.
\end{itemize}
\end{dfn}

A Morse function $f:\R^2\to\R$ is a function which has only
quadratic singularities. A {\it maximum} (resp. {\it saddle} and
{\it minimum}) of $f$ is a quadratic singularity with Morse index
$0$ (resp. $1$ and $2$). A {\it level set} (or an {\it $r$-level set})
of $f$ is the set in $\R^2$ given by $X_r:=\{x\in \R^2\mid f(x)=r\}$
for $r\in\R$. In particular, if $r$ is a regular value
then $X_r$ consists of disjoint smooth curves in $\R^2$ and
if $r$ is a critical value of only saddle singularities then
$X_r$ consists of immersed curves in $\R^2$.

\begin{dfn}
Let $P$ be an admissible divide in $D_n$ and fix
a checkerboard coloring of $D_n\setminus P$.
A {\it Morse function $f_P$ associated with $P$} is a Morse function
$f_P:\R^2\to\R$ which satisfies the following:
\begin{itemize}
   \item[(1)] the $0$-level set of $f_P$
coincides with $P$ in $D_n$ and, in particular,
$0\in\R$ is either a regular value or a critical value of
only saddle singularities of $f_P$;
   \item[(2)] each interior region of $P$ with black (resp. white) color
contains one maximum (resp. minimum) of $f_P$ and, in its small neighborhood
with coordinates $(x_1,x_2)$, $f_P$ is locally given by
$f_P(x_1,x_2)=x_1^2+x_2^2$ (resp. $f_P(x_1,x_2)=-x_1^2-x_2^2$);
   \item[(3)] each double point of $P$ corresponds to a saddle of $f_P$
and, in its small neighborhood with coordinates $(x_1,x_2)$,
$f_P$ is locally given by $f_P(x_1,x_2)=x_1^2-x_2^2$;
   \item[(4)] there are no singularities of $f_P$ in $D_n$
other than those in (2) and (3);
   \item[(5)] if the outside boundary of an exterior region of $P$
is exactly a component of $\bd D_n$ then the outside
boundary is contained in a level set of $f_P$;
   \item[(6)] if an exterior region of $P$ is not in case (5) then
there is just one point, in its outside boundary, at which
a level set of $f_P$ intersects $\bd D_n$ tangentially.
\end{itemize}
\end{dfn}

\begin{figure}[htbp]
   \centerline{\input{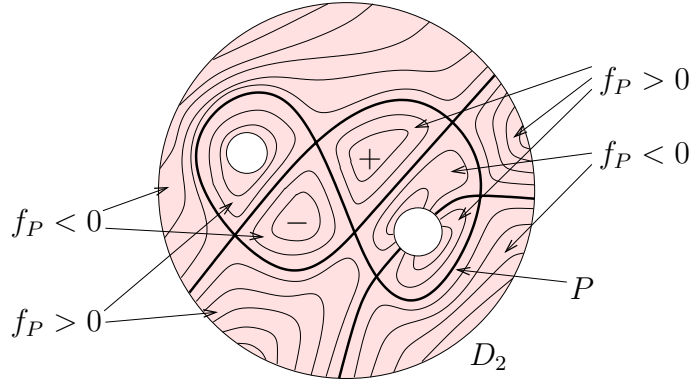}}
   \caption{An example of an admissible divide $P$ in $D_2$ with level sets
which determine a Morse function $f_P$ associated with $P$.\label{stein1}}
\end{figure}

For each admissible divide we can describe level sets as shown
in Figure~\ref{stein1}. The level sets determine a Morse function
$f_P$ of $P$ satisfying the above conditions.

For defining a complex valued function $F_P$ we prepare a few notations.
Let $a_i$, $i=1,\cdots\ell$, be the critical points of $f_P$.
Suppose that $U_i=\{x\in D_n\mid |x-a_i|<r_i\}$, for some $r_i>0$,
is the small neighborhood of $a_i$ introduced in the condition (2) or (3)
above, and set $U'_i=\{x\in D_n\mid |x-a_i|<r'_i\}$,
where $r_i>r'_i>0$. Let $\chi:D_n\to [0,1]$ be a positive
$C^\infty$-differentiable bump function such that $\chi(x)=0$ for
$x\in D_n\setminus\cup_{i=1}^\ell U_i$ and
$\chi(x)=1$ for $x\in\cup_{i=1}^\ell U'_i$.

The complex valued function $F_P$ is a map from $T(\R^2)$,
which is the tangent bundle to $\R^2$, to $\C$ defined by
\begin{equation}\label{fibmap}
   F_P(x,u):=f_P(x)+idf_P(x)(u)-\frac{1}{2}\chi(x)H_{f_P}(x)(u,u),
\end{equation}
where $x\in\R^2$, $T_x(\R^2)$ is the set of tangent vectors to $\R^2$ at $x$,
$u\in T_x(\R^2)$, $i=\sqrt{-1}$, $df_P(x)$ is the differential
and $H_{f_P}(x)$ is the Hessian of $f_P$ at $x\in\R^2$.
By setting $(x,u)$ to be the complex valued coordinates $x+iu$
in $\C^2$ as before, we can regard $T(\R^2)$ as $\C^2$ and $F_P$ as
a function from $\C^2$ to $\C$.

Recall that $M_n$ is the $3$-manifold $\#nS^1\times S^2$ embedded in $\C^2$.

\begin{dfn}\label{dfn27}
The {\it link} $L(P)$ of an admissible divide $P$ in $D_n$
is the set in $M_n$ defined by
\begin{equation}\label{link1}
   L(P):=M_n\cap F_P^{-1}(0).
\end{equation}
It can also be defined by
\begin{equation}\label{link2}
   L(P):=\{x+iu\in M_n\mid x\in P,\, u\in T_x(P)\},
\end{equation}
where $T_x(P)$ is the set of tangent vectors to $P$ at $x\in D_n$.
\end{dfn}

The coincidence of~\eqref{link1} and~\eqref{link2} is proved
in~\cite[Lemma~4.1]{i}. The expression~\eqref{link1} is important
for observing the Lefschetz fibration $F_P:N_\delta(D_n)\to\C$, while
the expression~\eqref{link2} will be used when we need to see
the position of the link $L(P)$ relative to the contact structure $\xi$.

\begin{thm}[\cite{acampo:99,i}]\label{fibthm}
There is a $3$-manifold $\hat M_n$,
obtained by an embedded isotopy of $M_n$ in $\C^2$ fixing the set $L(P)$,
such that the argument map $\pi:\hat M_n\setminus L(P)\to S^1$
defined by $\pi:=F_P/|F_P|$ is a locally trivial fibration over $S^1$
and its monodromy diffeomorphism is a product of positive Dehn twists.
In particular, it gives a positive open book decomposition of $M_n$
with $L(P)$ the binding.
\end{thm}

To prove our main theorem, we need an explicit description of
a fiber surface of the open book decomposition in Theorem~\ref{fibthm}.
In the rest of this section we describe it according to the observation
in~\cite{acampo:98a}.

\begin{figure}[htbp]
   \centerline{\input{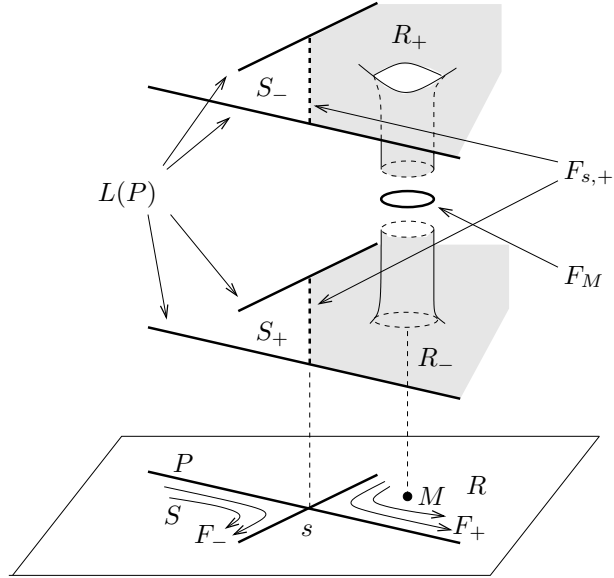}}
   \caption{Building up the fiber $F_1$.\label{stein2}}
\end{figure}

Let $P$ be an admissible divide in $D_n$.
The fiber surface over $1\in S^1$ is determined by the
equation $F_P/|F_P|=1$, which is equivalent to 
the condition that $F_P(x+iu)$ is real and positive.
We denote the fiber surface by $F_1$. Put
\[
   P_+:=\{x\in D_n\setminus\bd D_n\mid f_P(x)>0,\,df_P(x)\ne 0\}.
\]
The level sets of $f_P$ define oriented foliations $F_+$ and $F_-$
on $P_+$ such that $x+iu\in F_+$ (resp. $\in F_-$) if $df_P(x)(iu)>0$
(resp. $<0$). Put
\[
  P_{+,+}:=\{x+iu\in M_n\mid x\in P_+,\, u\in T(F_+)\}
\]
and
\[
  P_{+,-}:=\{x+iu\in M_n\mid x\in P_+,\, u\in T(F_-)\},
\]
where $T(F_\pm)$ is the set of tangent vectors on $P_+$ lying
in the same direction as $F_\pm$. Put
\[
   F_M:=\{x+iu\in M_n\mid x=M\}
\]
for each maximum $M$ in $P_+$, and also put
\[
   F_{s,+}:=\{x+iu\in M_n\mid x=s,\, H_{f_P}(x)(u,u)<0\}
\]
for each double point $s$ of $P$, which is a saddle point of $f_P$.
Finally put
\[
   \bd D_+:=\{x\in\bd D_n\mid f_P(x)>0\}.
\]
Then the fiber surface $F_1$ is given by the union
\[
   F_1=P_{+,+}\cup P_{+,-} \cup\bd D_+\cup
       \bigcup_{s\in P}F_{s,+}\cup\bigcup_{M\in P_+} F_M,
\]
where the gluings are ruled as shown in Figure~\ref{stein2}.
In the figure, $R$ is an interior region while $S$ may be an exterior region,
and $R_\pm$ (resp. $S_\pm$) is the lift of $R$ (resp. $S$)
corresponding to the foliation $F_\pm$.

\section{Oriented divides}

Before starting a proof of the main theorem,
we introduce oriented divides in $D_n$ and their links
defined in $\#nS^1\times S^2$. These will be used
for representing the given link $L$ in the main theorem
as we did in~\cite{gi} and~\cite{i}.

\begin{dfn}
An {\it oriented divide} $\vec P$ in $D_n$ is the image of a generic
immersion of a finite number of copies of the oriented unit circle into $D_n$.
Here the generic conditions are (i) and (iii) in Definition~\ref{dfndiv}.
\end{dfn}

\begin{dfn}
The {\it link} $L(\vec P)$ of an oriented divide $\vec P$ in $D_n$
is the set in $M_n$ defined by
\[
   L(\vec P):=\{x+iu\in M_n\mid x\in P,\, u\in T_x(\vec P)\},
\]
where $T_x(\vec P)$ is the set of tangent vectors to the immersed curve
of $\vec P$ at $x$ whose direction is consistent with the orientation
of $\vec P$ at $x$.
\end{dfn}

\section{Proof of Theorem~\ref{mainthm}}

The key point of the proof of main theorem is 
to apply the technique of Akbulut and Ozbagci in~\cite{ao:2001}
to the fiber surfaces of positive open book decompositions
of $\#nS^1\times S^2$ associated with divides in $D_n$.

We first state a theorem of Eliashberg which characterizes
compact Stein surfaces in words of $2$-handle attachings.
Let $W_n$ be a $4$-manifold obtained from $4$-ball by attaching
$n$ copies of $1$-handles.
Note that the boundary $\bd W_n$ of $W_n$ is homeomorphic
to $\#nS^1\times S^2$.
Let $K$ be a Legendrian knot in $\bd W_n$ with respect to
a contact structure.
The {\it canonical framing} of $K$
is the linking of the knot $K$ and a knot $K'$ obtained
by pushing-off $K$ into the direction normal to
the $2$-plane field of the contact structure.

The following is implicitly given by Eliashberg in~\cite{eliashberg:90}.
The statement can be found in~\cite{gompf:98}.

\begin{thm}[\cite{eliashberg:90}, cf.~\cite{gompf:98}]
An oriented, compact, smooth $4$-manifold with boundary is 
a Stein surface if and only if it can be obtained from
$W_n$, for some $n\geq 0$, by attaching $2$-handles along
Legendrian knots $K_i$, $i=1,\cdots,m$,
with respect to the standard contact structure
in $\#nS^1\times S^2$, with coefficient the canonical
framing minus $1$.
\end{thm}

Note that every Stein fillable $3$-manifold
can be obtained as the boundary of a compact Stein surface.

We need to introduce Lefschetz fibrations for relating
the above handle attachings to positive open book decompositions.
Let $W$ be a compact, connected, oriented, smooth $4$-manifold.
A {\it Lefschetz fibration} $p:W\to\Delta$ over a disk $\Delta$
is a map such that
\begin{itemize}
   \item each critical point $a$ of $p$ lies in the interior of $W$
and admits a coordinate neighborhood with complex valued coordinates
$z:=(z_1,z_2)$ consistent with the given orientation of $W$;
   \item in a small neighborhood of each critical point $a$,
the map $p$ is locally given by the form $p(z)=p(a)+z_1^2+z_2^2$;
   \item the orientation of $p(z)$ is consistent with that of $\Delta$;
   \item each fiber $p^{-1}(c)$, $c\in\Delta\setminus\bd\Delta$,
intersects $\bd W$ transversely;
   \item $p^{-1}(\bd\Delta)$ is contained in $\bd W$.
\end{itemize}
A Lefschetz fibration is called {\it allowable} if all its vanishing
cycles are homologically non-trivial in the fiber.
Note that a Lefschetz fibration $p:W\to\Delta$ induces,
on its boundary $\bd W$, a positive open book decomposition of $\bd W$
with binding the set $p^{-1}(0)\cap \bd W$.

The total space of a Lefschetz fibration has
a canonical handle decomposition characterized by
the critical values of the fibration map, see~\cite{kas:80,ao:2001}.
The following lemma gives a method for making
Lefschetz fibrations by using $2$-handle attachings.

\begin{lemma}[\cite{ao:2001} Remark~1]\label{lemma_ao}
Let $p:W\to\Delta$ be an allowable Lefschetz fibration over a disk $\Delta$.
Suppose that $W'$ is obtained from $W$ by attaching a $2$-handle
with coefficient the product framing minus $1$.
Then this $2$-handle attaching induces an allowable Lefschetz fibration
$p':W'\to\Delta'$ over a disk $\Delta'$.
\end{lemma}

In~\cite{ao:2001}, Akbulut and Ozbagci gave a proof of the existence
of a positive open book decomposition for every Stein fillable $3$-manifold.
The main idea in their proof is
to find a good position of the fiber surface of a torus knot in $S^3$
such that, for each $i=1,\cdots,m$,
$K_i$ stays on the surface and the product framing of $K_i$
coincides with $tb(K_i)$.
We will use the same strategy later.

Now we give a proof of the main theorem.

\vspace{5mm}
\noindent
{\it Proof of Theorem~\ref{mainthm}.}~~
Let $W$ be a Stein surface whose boundary is
the prescribed Stein fillable $3$-manifold, denoted by $M$.
Let $L$ denote the given link in $M$.
We suppose that $W$ is obtained from $W_n$
by attaching $2$-handles $H_i$ along Legendrian knots $K_i$,
with respect to the standard contact structure in $\bd W_n$,
with coefficients $tb(K_i)-1$, where $i=1,\cdots,m$.
Using regular isotopy, we can isotope the link $L$ in $\bd W$ so that
it does not intersect the $2$-handles $\cup_{i=1}^m H_i$.
Thus we can assume that $L$ is a link in
$(\#nS^1\times S^2)\setminus (\cup_{i=1}^m K_i)$.

Recall that the $3$-manifold $M_n$ is the boundary of a small compact
neighborhood of $D_n$ in $\C^2$ and homeomorphic to $\#nS^1\times S^2$.
We have seen that $M_n$ has the decomposition
$M_n=\bd N_A\cup \bd N_B$ and the piece $\bd N_B$ can be regarded as
the unit tangent bundle to $B$, where $B\subset D_n$ is the complement
of a small compact neighborhood of $\bd D_n$ in $D_n$.
As mentioned in Lemma~\ref{lemma1},
each Legendrian knot $K_i$ has a wave front representative $w_i$ on $B$.
We set $K:=\cup_{i=1}^m K_i$ and $w:=\cup_{i=1}^m w_i$.
We can also get an oriented divide representative of the link $L$ on $B$.

\begin{claim}\label{lemma2}
Every link in $M_n$ has an oriented divide representative in $B$
up to regular isotopy in $M_n\setminus K$.
\end{claim}

\begin{proof}
The claim is similar to Theorem~1.1 in~\cite{gi}
and can be proved by using the same technique (cf.~\cite{i,cgm}).
We here explain roughly how to prove this claim.
Let $L$ denote the given link in $M_n$.
Using regular isotopy, we assume that $L$ is included in $\bd N_B$.
If we think of $\bd N_B$ as the unit tangent bundle to $B$,
the link $L$ corresponds to a continuous assignment of
vectors with length $\delta>0$ and based on a curve $C$ in $B$.
By setting $L$ in a generic position, we can assume that the curve $C$
consists of generically immersed circles in $B$.
Assign an orientation to each component of $C$
and isotope the vectors so that they are in the same direction as
the orientation of $C$ except over a finite number of short intervals.
Here we first isotope the vectors in a neighborhood of
the double points of $C$
and in a neighborhood of the intersection points of $C$ and $w$
and then isotope the vectors on the other parts
so that the short intervals are included in
$C\setminus (\{double~points\}\cup w)$.
The proof is done by replacing each short interval by a small loop
such that the move of tangent vectors to the small loop
is consistent with the move of the vectors of
the isotoped continuous assignment of vectors corresponding to the link $L$.
\end{proof}

Thus we have a wave front representative $w$ of $K$
and an oriented divide representative, denoted by $\vec P$, of $L$ on $B$.
Using Legendrian isotopy for $K$ and regular isotopy for $L$,
we can assume that the union of immersed curves of these representatives
has only node and cusp singularities. Note that
the co-orientation of $w$ and the orientation of $\vec P$ can not
be in the same direction at each
intersection of the immersed curves of $w$ and $\vec P$
since the links $K$ and $L=L(\vec P)$ do not intersect each other.

Now we regard the immersed curve of $\vec P$ as a divide in $D_n$
and denote it by $P$.
We then add another divide $P'$ to $P$ so that $P\cup P'$ satisfies
the admissible conditions.
Let $a_1,\cdots,a_r$ be the intersection points of $w$ and $\vec P$,
$b_1,\cdots,b_s$ the self intersection points of $w$, and
let $c_1,\cdots,c_t$ be the cusp singularities of $w$.
Let $U_j$ (resp. $V_k$ and $W_\ell$) be a sufficiently small
neighborhood of $a_j$ (resp. $b_k$ and $c_\ell$).

\begin{figure}[htbp]
   \centerline{\input{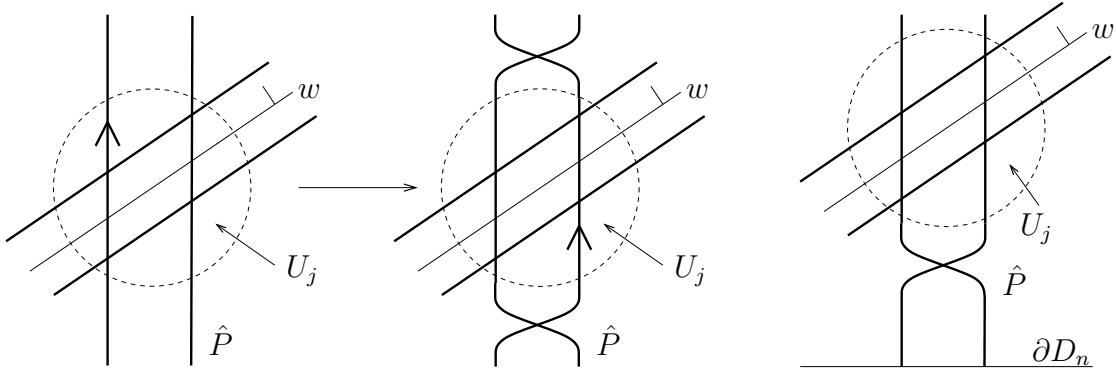}}
   \caption{The left and center figures show
the addition of two crossings in both sides of $U_j$.
The divide $\hat P$ is represented by thick curves
and $w$ is represented by co-oriented thin ones.
The orientation written on the curve of $\hat P$ represents
that of $\vec P$. The right figure shows the addition of a crossing between
$U_j$ and $\bd D_n$.\label{stein9}}
\end{figure}

\begin{figure}[htbp]
   \centerline{\input{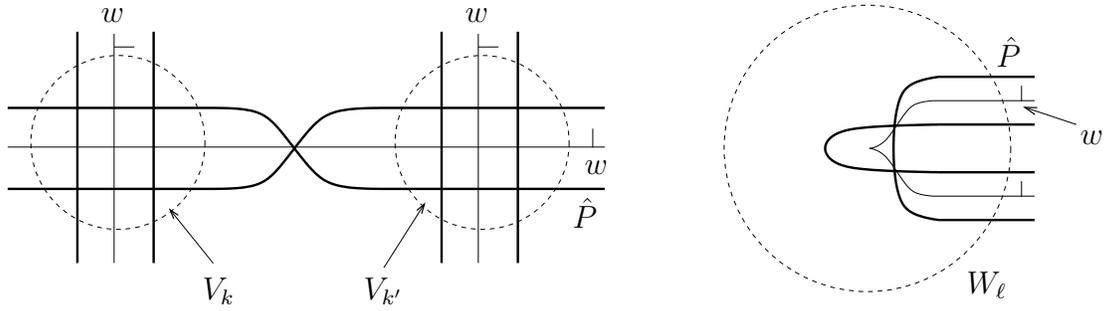}}
   \caption{The left figure shows the addition of a crossing to
the two parallel curves between $V_k$ and $V_{k'}$.
The right figure is the divide $\hat P$ inside $W_\ell$.\label{stein3}}
\end{figure}

Next we make a new divide $\hat P$ from $P\cup P'$ according to
the following rules:
\begin{itemize}
   \item[(1)] add an immersed curve parallel to $P\cup P'$;
   \item[(2)] outside $W_\ell$, write two immersed curves
parallel to $w$ in such a way that $w$ lies between the two parallel curves;
   \item[(3)] in $U_j$, if the argument between the orientation of $\vec P$
and the co-orientation of $w$ is less than $\pi/2$ and
$\vec P$ stays in the left-side line with respect to the co-orientation of $w$,
then add two crossings to the vertical two parallel curves near $U_j$,
as shown on the left and center in Figure~\ref{stein9},
so that $\vec P$ stays in the right-side line in $U_j$;
   \item[(4)] if a region of the divide obtained according to
steps (1) and (2) intersects a neighborhood $U_j$ and the boundary $\bd D_n$ then
add a crossing to the two parallel curves
as shown on the right in Figure~\ref{stein9};
   \item[(5)] if a region of the divide obtained intersects
two neighborhoods $V_k$ and $V_{k'}$ (possibly $V_k=V_{k'}$)
then add a crossing to the two parallel curves
as shown on the left in Figure~\ref{stein3};
   \item[(6)] inside each $W_\ell$, we define $\hat P$ as shown
on the right in Figure~\ref{stein3}.
\end{itemize}
Note that the divide $\hat P$ is admissible and contains $P\cup P'$.
In particular, the link $L=L(\vec P)$ is still contained
in the binding of the positive open book decomposition
associated with the divide $\hat P$.

We remark that $\hat P$ and $w$ may intersect orthogonally in $V_k$.
In this case, since the co-orientation of $w$ is tangent to $\hat P$
at the orthogonal intersection point, the link $K$ intersects a component
of $L(\hat P)$, but this is not a component of $L=L(\vec P)$
since $\vec P$ does not pass through the neighborhood $V_k$.

By changing the sign of $f_{\hat P}$ if necessary,
we assume that $f_{\hat P}>0$ in the regions
between the two parallel curves described in the above steps
outside $U_j$, $V_k$ and $W_\ell$.
Let $F_1$ denote the fiber surface, over $1\in S^1$,
of the positive open book decomposition
associated with the admissible divide $\hat P$.
As explained in Section~3, the surface $F_1$ can be understood
as the closure of the foliations $F_{\pm}$ on the regions $P_+$.

We will move the surface $F_1$ so that the link $K$
stays on it. During the move, $L=L(\vec P)$ can not intersect $K$
though the other components of $L(\hat P)$ can intersect it.
Since it is difficult to draw figures of the surface moved,
we describe the position of $K$ relative to the moved surface,
that is, we fix the divide $\hat P$ and its foliations $F_{\pm}$,
which represent the surface $F_1$, and
describe the move of $K$ by a family $\hat w_t$, $t\in [0,1]$,
of immersed curves with continuous assignment of vectors such that $\hat w_0=w$,
where the continuous assignment of vectors is the co-orientation
of $w$ when $t=0$ but it is not necessary to be orthogonal
to the immersed curves when $t>0$.
We call such a curve a {\it swinging wave front representative}.

Since $L=L(\vec P)$ can not intersect $K$ during the move,
the swinging co-orientation of $\hat w_t$ can not lie in the same direction
as the orientation of $\vec P$
at each intersection point of $\vec P$ and $\hat w_t$.
Also, at each self intersection point of $\hat w_t$
the two swinging co-orientations can not be in the same direction
otherwise the isotopy type of $K$ may change.

\begin{figure}[htbp]
   \centerline{\input{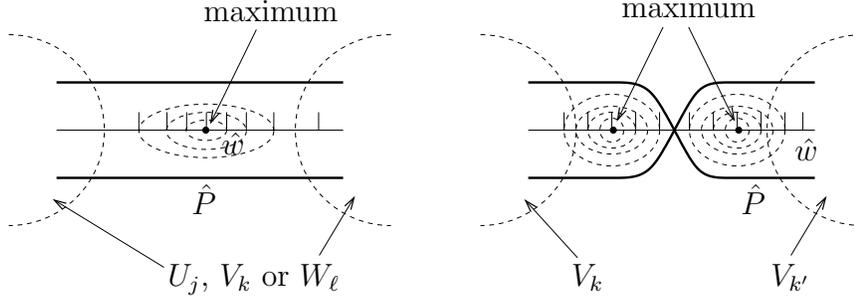}}
   \caption{The swinging wave front representative $\hat w$
outside $U_j$, $V_k$ and $W_\ell$.\label{stein4}}
\end{figure}

Now we move the surface $F_1$ to another surface $\hat F_1$
under the above rules so that $K$ is realized, at $t=1$,
by the swinging wave front representative $\hat w$
($=\hat w_1$) described in the following:
If there exists the curve of $w$ in an interior region not included
in $U_j$, $V_k$ and $W_\ell$ then we set $\hat w$
as shown in Figure~\ref{stein4}.
Here we do not need to swing the co-orientation, i.e.
we do not need to move the surface $F_1$.
The dots represent the maxima of the interior regions
and the dotted circles centered at the maxima are
the level sets of $f_{\hat P}$. The wave front representative $\hat w$
is, at each point, either tangent to a level set of $f_{\hat P}$ or
passing though a maximum. Hence the link $K$
is on the fiber surface $\hat F_1$ in this part.

\begin{figure}[htbp]
   \centerline{\input{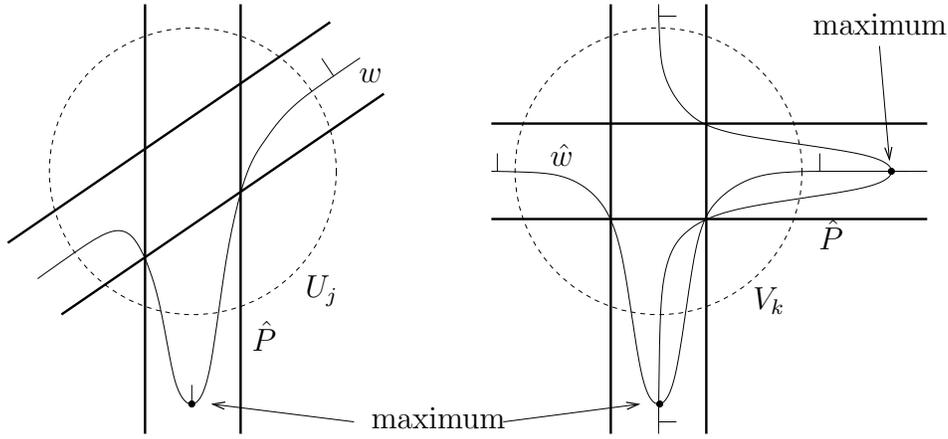}}
   \caption{The swinging wave front representative $\hat w$
inside $U_j$ and $V_k$.\label{stein5}}
\end{figure}

We set $\hat w$ inside $U_j$ as shown on the left in Figure~\ref{stein5}.
If $w$ runs from top-left to bottom-right in $U_j$ then we
need to consider the mirror image of the figure.
For making the figure to be simple, we did not draw
the level sets of $f_{\hat P}$, cf.~Figure~\ref{stein4}.
We set the swinging co-orientations of $\hat w$ to be
tangent to the level sets of $f_{\hat P}$ outside the quadratic singularities.
Then by observing $\hat w$ and $F_{\pm}$
we can conclude that $K$ lies on $\hat F_1$ in this part.

We need to check that there is no intersection of $K$ and $L$ during
the move $\hat w_t$. An intersection of $K$ with $L(\hat P)$ occurs
if the swinging co-orientation of $w_t$ and the tangent direction to $\hat P$
coincide at an intersection point of $\hat w_t$ and $\hat P$.
By observing the family $\hat w_t$ in $U_j$,
we can verify that the intersection occurs once over
the left vertical curve of $\hat P$ described on the left in Figure~\ref{stein5}.
Since we have set in step~(3) that the curve of $\vec P$ is
the right vertical line in the figure, we can conclude that
there is no intersection of $K$ and $L=L(\vec P)$ over $U_j$ during this move.

We also need to check if there really exists the maximum which
$\hat w$ passes and if there is no conflict with other neighborhoods.
First note that the region containing the maximum does not connect
to $V_k$ and $W_\ell$ since
$\hat w$ does not passes though $U_j$ vertically in the figure.
If the region is connected to the boundary $\bd D_n$ then
there is a maximum due to step~(4).
The rest is the case where it is connected to a neighborhood $U_{j'}$.
If $U_j=U_{j'}$ then the union of the region and $U_j$
constitutes an embedded annulus in $D_n$ and $\hat w$ intersects
it once transversely, which contradicts the fact that $\hat w$
consists of closed curves. Hence we can assume that $U_j\ne U_{j'}$.
If the co-orientation of $w$ in $U_{j'}$ is in the same direction
as that of $w$ in $U_j$ (along the vertical curves of $\hat P$ in the figure)
then $\hat w$ in $U_{j'}$ passes through a maximum different from
the one corresponding to $U_j$.
Hence we can observe $U_j$ and $U_{j'}$ independently.
If the co-orientation of $w$ in $U_{j'}$ is in the opposite direction
to that of $w$ in $U_j$ then $\hat w$ in $U_{j'}$ passes through
the same maximum as the one corresponding to $U_j$.
But we can also observe $U_j$ and $U_{j'}$ independently
since the co-orientations are in the opposite directions.

Inside $V_k$, we set $\hat w$ as shown on the right in Figure~\ref{stein5}.
The figure of $V_k$ is in the case where
the co-orientation of the vertical component of $w$ is in the right direction.
If it is opposite then we need to consider the mirror image of the figure.
We can see that $K$ lies on $\hat F_1$ as in the case of $U_j$
and that there is no intersection of $K$ and $L$ during the move $\hat w_t$
since $\vec P$ does not pass through a neighborhood of type $V_k$.
We can also check that the two strands of $K$ corresponding to
the vertical and horizontal curves of $\hat w$ in $V_k$
do not intersect each other during this move.
Since $w$ passes between both the vertical and horizontal parallel curves
of $\hat P$ in $V_k$,
these curves connect to either the horizontal curves in $U_j$
on the left in Figure~\ref{stein5}, or parallel curves in $V_{k'}$ (possibly
$V_k=V_{k'}$), or those in $W_\ell$.
Obviously, the maxima which $\hat w$ passes do not conflict with
those corresponding to $U_j$, and they also do not conflict with
$V_{k'}$ by the setting in step~(5).
Hence we can observe each $V_k$ independently.

\begin{figure}[htbp]
   \centerline{\input{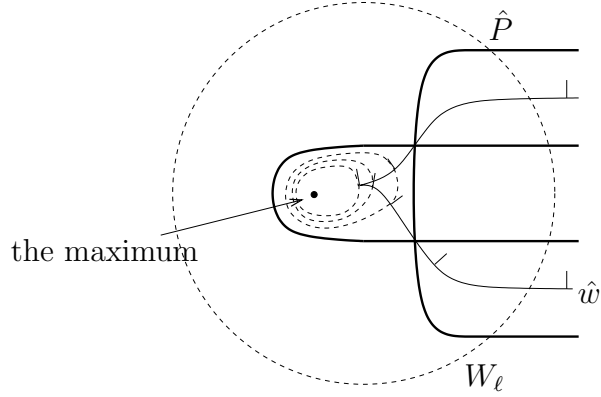}}
   \caption{The swinging wave front representative $\hat w$
inside $W_\ell$.\label{stein6}}
\end{figure}

Inside $W_\ell$, we set $\hat w$ as shown in Figure~\ref{stein6}.
The figure is in the case where the co-orientation is in the upper
direction. If it is opposite then we change the co-orientation
in the figure into the opposite direction.
As in the previous cases, we can easily check that $K$ lies on $\hat F_1$
in the figure.

Thus we got the move from $F_1$ to $\hat F_1$
without intersection of $K$ and $L=L(\vec P)$ such that $K$ lies on
the surface $\hat F_1$.

The product framing of $K_i$ associated with the Lefschetz fibration
is the framing of $K_i$ obtained by pushing-off it
in the direction normal to the fiber surface $\hat F_1$.
We now prove that this framing coincides with
the canonical framing of $K_i$.

\begin{claim}\label{claim54}
For each $i=1,\cdots,m$, the product framing of $K_i$ coincides with
the canonical framing of $K_i$.
\end{claim}

\begin{proof}
Let $v$ be a component of $\hat w$.
We assign an orientation to $v$ and put the labels
$p_1,q_1, p_2,q_2,\cdots,p_r,q_r$ to the maxima of $f_{\hat P}$ and
the cusps and double points of $\hat P$ according to the orientation.
Since the double points appear at every other label,
we can assume that $p_i$ is either a maximum or a cusp,
and $q_i$ is a double point.
We then separate $v$ into $r$ intervals $I_1,\cdots,I_r$ in such a way that
$I_i$ contains the pair of points $p_i$ and $q_i$.

We first prove the coincidence of the framings for each interval $I_i$
with $p_i$ a maximum. Instead of pushing-off the link of $v$ in
the direction normal to $\hat F_1$, we shift it on $\hat F_1$ into
one of the possible directions. This can be represented
by using a swinging wave front representative
as shown on the top in Figure~\ref{stein7}.
This swinging wave front representative can isotope as shown
on the bottom in the figure. By setting the isotoped swinging co-orientation
to be orthogonal to the immersed curve,
we can see that the curve is the one obtained by pushing-off
the swinging wave front representative $\hat w$
in the direction of the co-orientation,
which means that the corresponding link can be obtained
from the link of $\hat w$ by pushing-off it in the direction normal to
the $2$-plane field $\xi$. Hence the two framings coincide in this part.

\begin{figure}[htbp]
   \centerline{\input{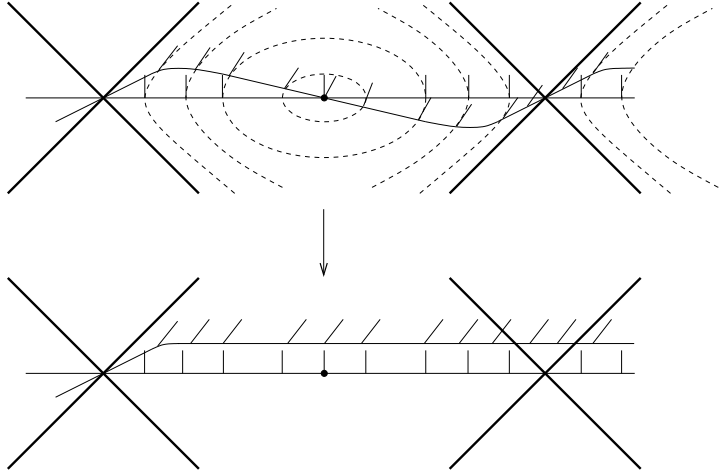}}
   \caption{Shift $\hat w$ on $\hat F_1$.\label{stein7}}
\end{figure}

Around each $V_k$ there are three self intersection points of $\hat w$
as seen on the right in Figure~\ref{stein5},
and there is no other self intersection of $\hat w$ in $D_n$.
For each $V_k$, we can consider the two curves of $\hat w$ independently
since the difference of the co-orientations of these two curves are
at least $\pi/2$ degree and we need to swing them a very little
during the shift in Figure~\ref{stein7}.
Hence the coincidence of the framings follows from the above observation.

The coincidence of the framings for each interval $I_i$
with $p_i$ a cusp can be proved by the same argument.
The shift of $v$ on $\hat F_1$ is represented
as shown in Figure~\ref{stein8}~(a).
\begin{figure}[htbp]
   \centerline{\input{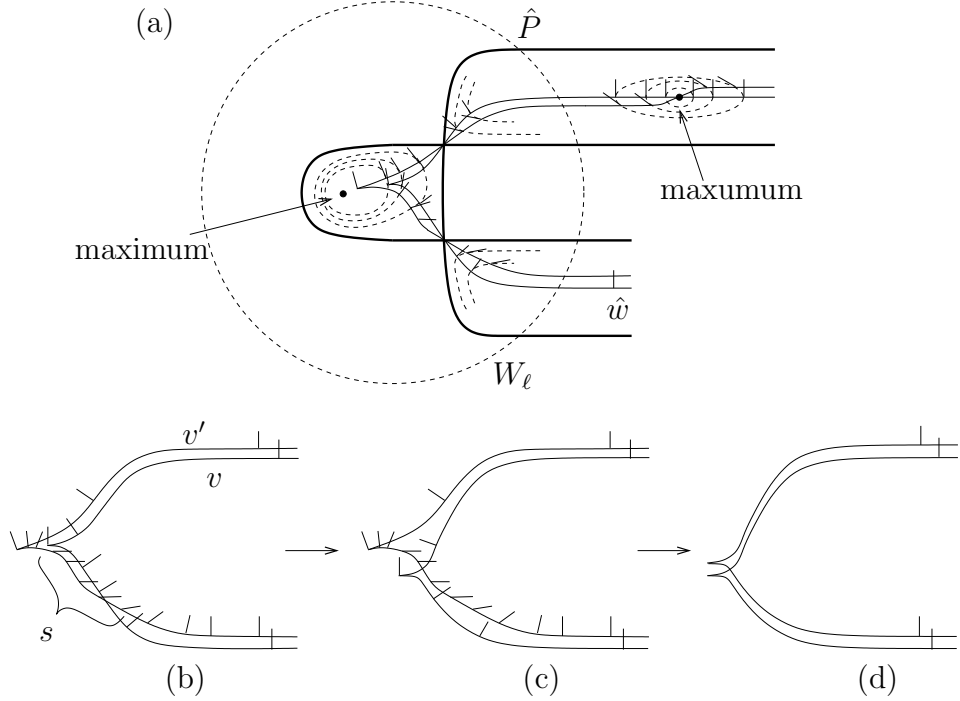}}
   \caption{Shift $\hat w$ on $\hat F_1$ near a cusp.\label{stein8}}
\end{figure}
The isotopy move of the swinging wave front representative from (a) to (b)
is the same move as in Figure~\ref{stein7}.
Denote by $v'$ the swinging wave front representative obtained
by the shift, see the figure~(b).
We assume that on the segment $s$ in $v'$ specified in the figure~(b)
the swinging co-orientation lies in the right direction.
Then we move the curve of $v$ so that the cusp comes under the cusp of $v'$.
This move does not change the linking of the links of $v$ and $v'$ since
all the swinging co-orientations on $s$ lie in the right direction
and no swinging co-orientation of $v$ lies in that direction.
Thus we obtain the figure~(c).
Finally we obtain the figure~(d) by swinging the co-orientations
of $v$ and $v'$ so that they are orthogonal to the immersed curves
of $v$ and $v'$ respectively.
This move does not change the linking of the links of $v$ and $v'$ since
the two co-orientations at the intersection point of $v$ and $v'$
do not intersect each other during the move.
Thus we can see that the curve is the one obtained by pushing-off
the swinging wave front representative $\hat w$
in the direction of the co-orientation, and hence the two framings coincide.
\end{proof}

We continue the proof of main theorem.
In~\cite{i} we introduced the Lefschetz fibration
associated with a divide in $D_n$ and obtained the fibration
in Theorem~\ref{fibthm} on its boundary,
and we can easily check that this Lefschetz fibration is allowable.
Hence, by the coincidence of the two framings in Claim~\ref{claim54},
we can conclude that the $2$-handle attachings along $w$, which yield
the compact Stein surface $W$, also produce a new allowable
Lefschetz fibration by Lemma~\ref{lemma_ao}.
The boundary of $W$ is the Stein fillable $3$-manifold $M$ expected
and it has a positive open book decomposition associated
with the Lefschetz fibration.
As already mentioned, the link $L$ is contained in the binding
of the positive open book decomposition of $M_n$ associated with $\hat P$.
Since the $2$-handle attachings do not change the binding,
$L$ is still contained in the binding of the positive open book
decomposition of $M$ associated with the Lefschetz fibration.
Suppose $L\cup L'$ is the binding. Then, as we did in~\cite{i},
we can modify the positive open book decomposition by plumbing positive Hopf
bands so that the binding is the union of $L$ and a knot in $M\setminus L$
(cf.~\cite{stallings:78,harer:82,gabai:83}).
This completes the proof of Theorem~\ref{mainthm}.
\qed

\section{A conjecture}

The main theorem of this paper suggests that we can freely choose the bindings
of positive open book decompositions of Stein fillable $3$-manifolds
and the level of freedom is higher than the complexity of links
in the $3$-manifolds.
Since every closed, oriented, smooth $3$-manifold has an
open book decomposition\cite{alexander:23},
it is natural to ask if any $3$-manifold has such freedom of the bindings.
We here propose a conjecture concerning this question.

\begin{conj}
Let $\mathcal M_k$, $k\geq 0$, be the set of closed,
oriented, smooth $3$-manifolds which have open book decompositions
whose monodromies contain at most $k$ negative Dehn twists.
Suppose that a $3$-manifold $M$ is in $\mathcal M_k\setminus\mathcal M_{k-1}$
(here we set $\mathcal M_{-1}=\emptyset$).
Then for any link $L$ in $M$ there exists a knot $L'$ such that
$L\cup L'$ is the binding of an open book decomposition of $M$
whose monodromy contains exactly $k$ negative Dehn twists.
\end{conj}

When $k=0$, $\mathcal M_0$ is the set of Stein fillable $3$-manifolds
and hence the conjecture is true by our main theorem.
However we know nothing about the class $\mathcal M_k\setminus\mathcal M_{k-1}$
for $k\geq 1$.
For instance, it is already an interesting problem to ask
the minimal numbers of negative Dehn twists
of open book decompositions of the non-Stein fillable $3$-manifolds
known in~\cite{lisca:98, lisca:99}.

\end{document}